\documentclass{aupl} 

\usepackage{
amsfonts,
latexsym,
amssymb,
enumerate,
verbatim,
mathrsfs,
color,
tikz,
}

\newcommand{\labbel}[1]{\label{#1} [[{\bf #1}]]}  
\newcommand{\bibbitem}[1]{\bibitem{#1} [[{\bf #1}]]}  
\renewcommand{\labbel}{\label} \renewcommand{\bibbitem}{\bibitem}  

\newcommand\circled[1]{\mathbf{#1}}

\usepackage{tikz} \renewcommand*\circled[1]{\tikz[baseline=(char.base)]{    \node[shape=circle, draw, inner sep=1pt,         minimum height=12pt] (char) {#1};}}

\newcommand\dad{\sum^{\Delta}}
\newcommand\daf{\sum^{\Phi}}
\newcommand{\ceq }{\mathrel{\stackrel{\scriptstyle  C}{=}}}

\numberwithin{equation}{section}

\newtheorem{theorem}{Theorem}[section]
\newtheorem{lemma}[theorem]{Lemma}

\newtheorem{proposition}[theorem]{Proposition}

\newtheorem{fact}[theorem]{Fact} 
 
\newtheorem*{claim*}{Claim}

\newtheorem*{theorem*}{Theorem}
\newtheorem*{proposition*}{Proposition}
\newtheorem*{corollary*}{Corollary}
\newtheorem*{lemma*}{Lemma}
\newtheorem*{scholion*}{Scholion}

\theoremstyle{definition}
\newtheorem{definition}[theorem]{Definition}

\newtheorem{conjecture}[theorem]{Conjecture} 
 
\newtheorem{problem}[theorem]{Problem}

\theoremstyle{remark}
\newtheorem{remark}[theorem]{Remark}

\newtheorem*{remark*}{Remark}
\newtheorem*{remarks*}{Remarks}
 
\newtheorem{example}[theorem]{Example}

\newtheorem*{observation*}{Observation}

 \allowdisplaybreaks[1]

\begin{document}

 \title{Infinite sums of combinatorial games (Dadaist games)}
 
\author{Paolo Lipparini} 
\urladdr{http://www.mat.uniroma2.it/~ lipparin}
\address{Dipartimento Giocoso di Matematica\\Viale della  Ricerca
Scientifica\\Universit\`a di Roma ``Tor Vergata'' 
\\I-00133 ROME ITALY
\\(currently retired)}


\thanks{Work performed under the auspices of G.N.S.A.G.A.}

\begin{abstract}
We propose an interpretation of the infinite sum
of a sequence of combinatorial games. In such an interpretation,
plays involve infinite runs, but without loops.
The notion of a run is quite natural, but different
possibilities arises for the notion of an
\emph{alternating run}. 
\end{abstract} 

\maketitle

\section{Introduction} \labbel{intro}

Conway \cite{onag,ww}
introduced a general notion of ``game'',
including, in particular,  surreal numbers,  
which generalize, at the same time,
real numbers and ordinals. Notions of
sums and limits for surreals numbers
have occasionally been studied, e.~g.\ 
 \cite[p. 40]{onag},
 \cite[Definition 3.19]{S}, \cite{LM,W}.
 See eg. \cite{A,Bac,Ch,Sier} for the special case of ordinal numbers and  \cite{Al,BEK,Be,Go,RS} for surreal numbers.
On the other hand, to the best of our knowledge, no
notion of limit or infinite sum has ever been proposed for arbitrary 
combinatorial games, before our preprint \cite{scg}.
 We show that, quite surprisingly, there is a natural notion of
an infinite sum of arbitrary combinatorial games, though
we need to allow infinite runs (but still no loop).

While we still know no special application
of this general notion of an infinite sum, it seems quite interesting
for itself and probably deserves further study.
The theory produces some quite counterintuitive 
conclusions; for example, an infinite sum of Left winner games
might turn out to be a second winner game.
This suggests that ``Dadaist games'' could be a possible name for such infinite sums of games.

Many problems are still open; in particular, 
while the notion of a (possibly infinite) run is quite natural,
there are many possible variations on the notion
of an alternating run. Moreover, in general,
we do not know whether
there is always a winning strategy for one of the two  players in 
a countable sum of games. The problem whether
some form  of  ``generalized'' Conway equivalence 
holds is still completely open.

We assume familiarity with the basic notions of 
combinatorial game theory.
Here a \emph{combinatorial game} is a position in a possibly infinite
two-player game with perfect information, no chance element, no draw
and no infinite run, as studied, for example, in \cite[Chapter VIII]{S}.
Games are possibly \emph{partizan}
and we always assume \emph{normal} (not mis\'ere) \emph{play}.  
The two players are called \emph{Left}
and \emph{Right} and in each play on some game
they can move either as the first or as
the second player. We sometimes capitalize
\emph{First} and \emph{Second} in order
to denote the player who plays first (or second) 
in some specific play. 
 
Here games are considered up to isomorphism.
Isomorphism of combinatorial games is intended in an extensional sense:
$2= \{ \, 1 \mid \, \} $ is isomorphic to 
$1+1$ while, formally, in $2$ Left has exactly one option,
and in  $1+1$ she has two (isomorphic) options:
playing on the first $1$ or in the second $1$.  
Isomorphism of combinatorial games
is denoted by  $ \cong $.
Isomorphic combinatorial games will be considered here as the
same game; equivalence classes of combinatorial  games under isomorphism
are called \emph{forms} of games. 
On the other hand, the equivalence class of a combinatorial game
under Conway equivalence is called its \emph{value}.
Here we assume the notation from \cite{S},
except that Conway equivalence is denoted by
$\ceq $.

\section{Preliminaries} \labbel{prel} 

The reader familiar with the basic notions of 
combinatorial game theory as presented in
 \cite[Chapter VIII]{S} might skip the present section.
As mentioned, here combinatorial games are possibly infinite,
possibly \emph{partizan}
and we always assume \emph{normal} (not mis\'ere) \emph{play}.  

In detail, a \emph{combinatorial game} is a position in a
two-player game with perfect information, no chance element, no draw
and no infinite run.
In more detail, a combinatorial game is an ordered pair
consisting of the  \emph{Left options}  and 
of the \emph{Right options}. Left and Right
are the names of the two players and \emph{options}
are themselves games; the word option is just meant
to suggest that they are possible moves for some player.
In fact, here \emph{game}, \emph{position} and 
\emph{option} are formally synonyms. 
Generally, games are not assumed to be \emph{impartial}:
they are \emph{partizan} in the sense that the
sets of moves allowed to the two players are possibly distinct.  

A \emph{move} in a game $G$ is the choice of some option;
it is a \emph{Left move} if it is a Left option, and similarly for the Right.
The simplest game is $0$, in which  neither Left nor Right has a move.
All the other games can be constructed, up to extensional isomorphism,
by a transfinite induction. For example, $1=\{ \, 0 \mid  \,\}$ 
is the game in which Left has the only option $0$ and Right cannot move.
Similarly,  in $-1=\{ \, \mid 0 \,\}$ only Right has a legal move.
In $*=\{ \, 0 \mid 0 \,\}$ both players have the only option $0$.
The games $0$ and $*$ are impartial, while $1$ and $-1$ are 
strictly partizan.     
 
Given a game $G$, a \emph{finite run}  is a sequence $G_0, G_1, G_2, \dots, G_n$
of games such that $G=G_0$ and each $G_{i+1}$ is an option of
$G_i$.  All the games $G_n$ which can be reached in this way, starting
from $G$, are called \emph{subpositions} of $G$, with $G$ 
being an \emph{improper} subposition of itself.
Combinatorial games are possibly infinite, but are assumed to have
no infinite run. An \emph{alternating} run is a run in which
Left and Right play alternatively. A \emph{play (on $G$)} is an alternating run
in which the player to move cannot move on 
  the last element $G_n$ of the run. In such a situation, 
the game is won by the Other player. Namely, a player
who cannot move at her turn loses.

A straightforward induction shows that combinatorial
games divide in four nonoverlapping categories:
the \emph{First winner}  games in which the player who starts playing
has a winning strategy, no matter whether First is 
Left or Right; the \emph{second winner}   games in which the second 
player has a winning strategy; the \emph{Left winner}  games in which
Left has a winning strategy, no matter  who begins
and, finally, the \emph{Right winner}  games. Thus each 
combinatorial game $G$ has
exactly one \emph{outcome} $o(G)$ among $ \mathcal N$ (first, the Next player),
$\mathcal P$ (second, the Previous player), $\mathcal L$ (Left)
and $\mathcal R$ (Right).  For example, $o(0)= \mathcal P$, 
$o(*)= \mathcal N$, $o(1)= \mathcal L$ and $o(-1)= \mathcal R$.

The main point is that, unlike classical games like chess, in which
the rules fix a first player, in the theory of combinatorial
games each player can be considered either as First, or as Second,
giving rise to two distinct plays, of course. This is necessary in order to
introduce the  sum of two games. 
The \emph{(disjunctive) sum} 
$G+H$ of the games $G$ and $H$ is the game
in which each player, at her turn, chooses exactly one among 
$G$ and $H$, moves on it, and the resulting position is, say,
$G^L+H$, if Left was to move and she has chosen the option
$G^L$ of $G$. Again, transfinite induction shows that
this is a good definition.

Isomorphism of combinatorial games,
denoted by  $ \cong $,  is defined recursively:
two games are \emph{isomorphic} if they
have the same \emph{sets} (i.e., multiplicities
do not count) of Left and Right options. 
Thus isomorphism is intended in an extensional sense:
as already recalled,
$2= \{ \, 1 \mid \, \} $ is isomorphic to 
$1+1$ while, formally, in $2$ Left has exactly one option,
and in  $1+1$ she has two (isomorphic) options:
playing on the first $1$ or in the second $1$.  
Isomorphic combinatorial  games will be considered here as the
same game; equivalence classes of games under isomorphism
are called \emph{forms} of games. 

As a much coarser relation, 
two combinatorial  games $G$ and $H$ are \emph{Conway equivalent},
written as $G\ceq H$,  
if $o(G+K) = o(H+K) $, for every game $K$.
The equivalence classes of combinatorial games
under Conway equivalence is called the class of \emph{values}
of combinatorial games and form a Group, with sum of games
as addition. 
If $G$ is a game, $-G$ is the game obtained from  
$G$ by recursively exchanging the Left and Right options
in every subposition of $G$.
Conway proved that $G\ceq  0$ if and only if
$o(G) = \mathcal P$ and that $G\ceq  H$
if and only if $G-H\ceq  0$.    

We refer to, e.~g., \cite{S} for full details about the above notions.

\section{Countable sums of combinatorial games} \labbel{count}

The definition of the (disjunctive) sum of two 
combinatorial games can be evidently
iterated a finite number of times.

Similar definitions can be applied to an infinite number of
combinatorial games. Informally,  any sequence
$( G_i) _{i \in I} $ of combinatorial games
can be considered as a ``$\Delta$-game'' in which   
each player, at her turn, chooses
some index $j \in I$ and picks  an option, say,
$G_j^L$ on $G_j$; the result of the move
is the sequence $( H_i) _{i \in I} $
in which all the games are unchanged, except for
$G_j$, which is replaced by  $G_j^L$.
As we will point out explicitly in 
Remark \ref{manydiff}(1), this definition cannot be 
given inductively, hence, formally, an infinite sum
of games has to be considered just as a sequence of games.
Since,  under the above intuitive notion,
infinite runs obviously may arise, we need also to deal explicitly with
the ``resultant''   of an infinite run, i.e., the
 ``remainder of the games after the run''.
See Definition \ref{run} for a precise definition. 

Such $\Delta$-games, as outlined above,
are not really well behaved.
For example, in the classical theory, second winner 
combinatorial games are neutral elements
for the Conway sum with respect to Conway equivalence.
This means that the outcome of some combinatorial game is not affected
if we add a second winner game.
On the other hand, as far as countable sums
are concerned, we will see that there is a  
second winner $\Delta$-game $\mathbf G$ such that
$\mathbf G + \mathbf H$ is still second winner,
\emph{for every} $\Delta$-game $\mathbf H$
(there is a natural definition of the sum of  $\Delta$-games,
see Definition \ref{sgam}(B)). 

We now present the formal definitions. 

\begin{definition} \labbel{sgam} 
(A) Let $( G_i) _{i \in I} $ be a sequence of combinatorial games,
from now on denoted as $\dad _{i \in I} G_{i}$,
and considered as the $\dad$-\emph{sum} of the  $G_i$'s.
Each $G_i$ will be called a \emph{summand} 
of $\dad _{i \in I} G_{i}$. The \emph{Left options}
of $\dad _{i \in I} G_{i}$
 are the sequences
$\dad _{i \in I} H_{i}$, with the same index $I$,
where $H_j$ is a Left option of $G_j$,
for some $j \in I$ and $H_i=G_i$,
for all $i \neq j$. \emph{Right options} of     
$\dad _{i \in I} G_{i}$ are defined symmetrically.
If no $G_i$ has  Left options, then
 $\dad _{i \in I} G_{i}$ has no Left option; similarly
for the Right options.

Sums of combinatorial games, as in the
above definition, will be called
\emph{$\Delta$-games}
and will be denoted by boldface upper case letters 
as $\mathbf  G$, $\mathbf H$, \dots\  
 For ease of notation, 
we frequently assume that $I = \mathbb N$, but notice
that the ordering on $\mathbb N$, or on any other index
set $I$, will usually play no role  
in the theory of $\Delta$-games.
So here is a difference in comparison with \cite{scg}.
 With the above convention,
we may write $G_0+G_1+G_2+ \dots +G_i+ \dots $
 or simply $G_0+G_1+G_2+ \dots $ 
for  $\dad _{i \in \mathbb N} G_{i}$.
We do not introduce a new symbol for the sum since,
 in the case when $I$ is finite,
the sum reduces to the usual disjunctive sum of combinatorial games.

(B)
There is an immediate way to define the sum 
$\daf$ 
of a family of $\Delta$-games.
If  $( \mathbf G_j) _{j \in J} $ is a sequence of $\Delta$-games,
with   $\mathbf G_j = \dad _{i \in I_j} G_i$, for each $j \in J$, then, assuming
without loss of generality that the sets $I_j$ are pairwise disjoint,
we let   $\daf _{j \in J}  \mathbf G_j = \dad _{i \in I} G_i$,
where $I=\bigcup _{j \in J} I_j $.  
In case  the sets $I_j$ are not pairwise disjoint, simply relabel
the indexes.
\end{definition}

Before continuing with formal definitions, 
we make some more or less informal comments.

\begin{remark} \labbel{manydiff}   
There are many important differences with respect to
the case of finite sums.

(1) As we mentioned in Section \ref{prel},  
(Forms of) combinatorial games can be 
constructed by a transfinite induction starting form
the zero game $0$. Hence every combinatorial game
has a rank, generally called its \emph{birthday}:
the first day (= ordinal) in which it has been created.   
In a sense, the theory of combinatorial games is completely self-contained;
it needs no mathematical foundational background, as far as
we are allowed to combine ``already existing games''. A formal axiomatization
is generally not explicitly presented, but it would be quite similar to axioms for
set theories, in the sense that games are some kind of ``bilateral'' or
``two-headed'' sets. 

As a basic illustration, for $G$ and $H$ 
combinatorial games, $G^L+H$ has birthday strictly smaller 
than $G+H$, hence a transfinite induction 
on birthdays can be generally carried over, in order
to prove theorems about finite sums of combinatorial games.

 On the other hand, here our definition of an infinite sum
of combinatorial games does rely on (non-game-theoretical) set theory:
the $\Delta$-game $\dad _{i \in I} G_{i}$
is defined as a sequence of combinatorial games, namely a function
from some set $I$ to the class of combinatorial games.
This is not just formal  idiosyncrasy:  as we will see below,
in the case of infinite
sums a $\Delta$-game might be isomorphic to some  subposition
of itself, hence $\Delta$-games cannot be constructed inductively.  
As a much more troublesome consequence, theorems about 
$\Delta$-games cannot be proved by transfinite induction.

(2) Obviously, $\Delta$-games generally have infinite runs.
For example,  let $ \mathbf  G$ be $1-1+1-1+1-1+\dots$,
an abbreviation for $1+(-1)+1+(-1)+1+(-1)+\dots$. 
A possible play on   $\mathbf G$, with Left playing as First, consists of the following
run:
\begin{equation}\labbel{1}
\begin{aligned} 
&1-1+1-1+1-1+\dots,
\\[-3pt]
&0-1+1-1+1-1+\dots,
\\[-3pt]
&0+0+1-1+1-1+\dots,
\\[-3pt]
&0+0+0-1+1-1+\dots,
\\[-3pt]
&0+0+0+0+1-1+\dots, 
&\qquad \dots
 \end{aligned} 
   \end{equation}    

But more involved runs are possible, for example,
\begin{equation}\labbel{2}
\begin{aligned} 
&1-1+ \underline{1}-1+1-1+1-1+1-1+1-1+\dots
\\[-3pt]
&1-1+0 -\underline{1} +1-1+1-1+1-1+1-1+\dots
\\[-3pt]
&1-1+0+0+1-1+\underline{1}-1+1-1+1-1+\dots
\\[-3pt]
&1-1+0+0+1-1+0-\underline{1}+1-1+1-1+\dots
\\[-3pt]
&1-1+0+0+1-1+0+0+1-1+\underline{1}-1+\dots
\\[-3pt]
&1-1+0+0+1-1+0+0+1-1+0-\underline{1}+\dots, 
&\qquad \dots,
 \end{aligned} 
 \end{equation}    
where we have underlined those summands on which a move is
done.

Can the  run \eqref{2}  be considered as a (complete) play? 
It is natural to assume that it is  incomplete;
indeed, it is natural to assume that after step $ \omega$ Left 
can still play on the first $1$ in the sequence, so that the play must go on. 
Thus we have to be precise about what ``happens''
after an infinite sequence of  moves and we have to define what is the 
 ``resultant'' of the run. This can be done in a rather
straightforward fashion, and it can be checked that $\Delta$-games actually
always stop after some run of  ordinal length, with a non ambiguous
(but not necessarily unique, see Sections \ref{pfs} and  \ref{reg})
possibility of declaring a ``winner''.  However, we cannot
perform proofs by induction on birthdays as in (1) above\footnote{We can,
however, perform proofs by induction on the length of a run.};
 in particular, we do not know
whether some Player has always  a winning strategy; though we will show
that a winning strategy exists in many cases.

(3) Consider now the $\Delta$-game
$ \mathbf  G = 1-1+0+1-1+0+1-1+0+\dots$
and the following run, where we group summands
three by three for ease of reading. 
\begin{equation}\labbel{3}
\begin{aligned} 
\ \ &1{-}1{+}0 \ {+}\  \underline{1}{-}1{+}0 \ {+}\  1{-}1{+}0 \ {+}\ 
1{-}1{+}0\ {+}\ 1{-}1{+}0\ {+}\ 1{-}1{+}0\ {+}\dots
\\
&1{-}1{+}0\ {+}\ 0{-}\underline{1}{+}0 \ {+}\   1{-}1{+}0 \ {+}\ 
1{-}1{+}0\ {+}\ 1{-}1{+}0\ {+}\ 1{-}1{+}0+\dots
\\
&1{-}1{+}0\ {+}\ 0{+}0{+}0 \ {+}\  1{-}1{+}0  \ {+}\ 
\underline{1}{-}1{+}0 \ {+}\ 1{-}1{+}0\ {+}\ 1{-}1{+}0+\dots
\\
&1{-}1{+}0\ {+}\ 0{+}0{+}0 \ {+}\   1{-}1{+}0\ {+}\ 
0{-}\underline{1}{+}0 \ {+}\ 1{-}1{+}0\ {+}\ 1{-}1{+}0+\dots
\\
&1{-}1{+}0\ {+}\ 0{+}0{+}0 \ {+}\   1{-}1{+}0\ {+}\ 
0{+}0{+}0 \ {+}\ 1{-}1{+}0\ {+}\ \underline{1}{-}1{+}0+\dots
\\
&1{-}1{+}0\ {+}\ 0{+}0{+}0 \ {+}\   1{-}1{+}0\ {+}\ 
0{+}0{+}0 \ {+}\ 1{-}1{+}0\ {+}\ 0{-}\underline{1}{+}0+\dots
\\
&\qquad \dots, 
\end{aligned} 
 \end{equation}    
 whose ``resultant'' intuitively is the 
following $\Delta$-game $\mathbf H$
\begin{equation*}
\begin{aligned}  
&1{-}1{+}0\ {+}\ 0{+}0{+}0 \ {+}\   1{-}1{+}0\ {+}\ 
0{+}0{+}0 \ {+}\ 1{-}1{+}0\ {+}\ 0{+}0{+}0+\dots,
\end{aligned} 
  \end{equation*}     
 where, again, we have underlined those summands on which a move is
done. Again intuitively, and full details will be presented soon,
the resultant $\mathbf H$ of the above run is isomorphic to $\mathbf G$.
This shows that we need to deal with \emph{sequences}
of games (thus indexed by some set), rather than with 
\emph{infinite sets} of games.
In the latter case, we could not see any difference
between the $\Delta$-games $\mathbf G$ and $\mathbf H$.
On the other hand, since we are keeping the index set fixed,
$\mathbf H$ is a proper subgame of $\mathbf G$. Thus in the theory  
of $\Delta$-games a game may be isomorphic to  some proper
subposition. But, again since we have an index set, every run
terminates after some set of moves, indexed by  a possibly infinite ordinal.
The reason is that each summand is a combinatorial game, hence, on each summand,
the play ends, sooner or later.

Hence the situation here is different with respect to 
the theory of \emph{loopy} games \cite[Chapter VI]{S}.
Technically, here the tree of the runs on some $\Delta$-game  
is well-founded.
 \end{remark}

Now for the full details.

\begin{definition} \labbel{run}   
(a) As in the classical treatment, a \emph{move} 
in a $\Delta$-game $\mathbf G=\dad _{i \in \mathbb N} G_{i}$
 is the choice of some option,
as defined in Definition \ref{sgam}. A move
 is a \emph{Left move} if it is a Left option, and similarly for the Right.
 
(b) For $\Delta$-games, the rigorous definition of a run
has to be given by a simultaneous ordinal induction
on the length of the run,
together with the corresponding definition of the
\emph{resultant} of the run.
The point is significant only  at limit stages. 
 
Given a $\Delta$-game $\mathbf G$, a run
on $\mathbf G$ 
  is a particular ordinal-indexed sequence
of $\Delta$-games; each of these games  will be called 
a \emph{subposition}  of $\mathbf G$.
The formal definition of a \emph{run on $\mathbf G$}
is given inductively as follows. 
Recall that an ordinal $\alpha$ is the set of all smaller
ordinals, hence the domain of an $\alpha$-indexed sequence
is the set of all ordinals which are \emph{strictly} smaller than $\alpha$.  

(Base case) A $0$-run $r$ on $\mathbf G$  is the empty sequence, the sequence
without elements;
 the resultant  of $r$ is $\mathbf G$ itself.

(Successor step)
 Suppose that  $\alpha$ is an ordinal and we already know what an
$\alpha$-run and its resultant are. An $\alpha+1$-run $s$ on $\mathbf G$ is an
$\alpha+1$-indexed sequence of 
$\Delta$-games
such that the restriction $r$  of 
$s$ to $\alpha$ is an $\alpha$-run on $\mathbf G$ and
$s( \alpha)$ is an  option $H$  in the 
 the resultant  of $r$. Here option is intended in 
the sense of Definition \ref{sgam}. The resultant
of $s$ is such an $H$. The move 
$s( \alpha)$ is called \emph{the $\alpha+1$th move in the run $s$.} 
Notice that, say, $0$ is the \emph{first} ordinal,
$1$   is the \emph{second} ordinal, etc.
This case is not essentially different from the definition
of a run in the standard case of combinatorial games.

(Limit step)
Suppose that $\alpha$ is a limit ordinal and we already know what a
$\beta$-run is, for $\beta < \alpha $, and we also
know what the corresponding resultants are. An  $\alpha$-run $s$ 
on $\mathbf G=\dad _{i \in \mathbb N} G_{i}$ is an
$\alpha$-indexed sequence of 
$\Delta$-games
such that, for every $\beta < \alpha $,  the restriction $r_ \beta $  of 
$s$ to $ \beta $ is a $ \beta $-run on $\mathbf G$.
We need also to define the resultant of $s$.
For every index $i$ and every $\beta < \alpha $,
the $i$th summand in the resultant  of $r_ \beta $
is a possibly improper subposition of $G_i$.
More generally, if $ \gamma  < \beta < \alpha $, then
the $i$th summand in the resultant of $r_ \gamma  $
is a subposition of  the $i$th summand in the resultant of $r_ \beta $
(formally, this should be proved simultaneously
with the present inductive definition, but we hope that this fact is evident).
Since all the summands in $\mathbf G$ are combinatorial games,
they have no infinite run; hence, since $\alpha$ is a limit ordinal, 
we have that, for every index $i$, the  $i$th summand in the resultants
is eventually constant, say, $G_i^*$. 
Then the resultant of $s$ is defined to be
$\dad _{i \in \mathbb N} G_{i}^*$.

If $\alpha$ is a limit ordinal, there is no such thing as the 
$\alpha$th move, just as there is no $0$th move.

(c) Having defined what a run is, we need some care
in defining alternating runs.
Recall that any ordinal $\alpha$ can be written
uniquely as $ \beta + n$, for $n \in \mathbb N$ and
$\beta$ either $0$ or a limit ordinal.  
 An \emph{alternating run 
in which Left plays as First}
is a run $s$ in which,
for every limit  ordinal $\beta$
(or  $ \beta = 0$), 
the $ \beta +n$ th  move
$s( \beta +n-1)$ 
 is a Left move,
for every odd $n$ and a Right move, for every 
even $n>0$. 
An alternating run 
in which Right plays as First is defined symmetrically.
An \emph{alternating run}
is a run in which either Left or Right plays as First.
Of course, in place of saying that ``Left  plays as First'',
we might say that ``Right  plays as Second''.

(d) We remark that, under the above definition, if, say, Left 
plays as First,
then, after a run of length $ \omega$, Left is still
the first to move. Similarly, after every run 
of limit ordinal length, the First Player will 
always be the first player to move.
This rule might be changed, possibly letting the roles of
First and Second alternate, or in various distinct fashions.
Some possibilities will be considered in the 
following sections.

(e) Finally, a \emph{play} is a maximal alternating run $r$,
that is, an alternating run in  which the player to move cannot move on 
  the resultant of $r$. As in the classical case, in such a situation, 
the game is won by the Other player.

The, possibly undefined,
 \emph{outcome} $o(\mathbf G)$  of a $\Delta$-game $\mathbf G$ 
is one among $ \mathcal N$ (first, the Next player),
$\mathcal P$ (second, the Previous player), $\mathcal L$ (Left)
and $\mathcal R$ (Right), as defined in the classical case of combinatorial games.   
Since here we cannot perform inductions on the complexity
of $\Delta$-games, it might happen that  there are games in which
neither player has a winning strategy. 
It is an open problem whether there are $\Delta$-games
with an undefined outcome.

(f) Suppose that $\mathbf G=\dad _{i \in  I} G_{i}$ 
and no $G_i$ is isomorphic to $0$.
If there is some run $r$ whose resultant is a sum
of games isomorphic to $0$, then $r$ induces a surjective
function $t$ from some ordinal to  $I$: just take 
 $t( \alpha )$ to be the index on which the move 
$r( \alpha )$ is done. 
If we take $G_i= \{ \, -1 \mid 1 \, \} $,
for every $i \in I$, then the existence of a play,
or just of some maximal run, implies that
$I$ is well-orderable.   
Thus the assumption that, 
for every $\Delta$-game, there exists a maximal run
for that game,  implies
that every set is well-orderable. Hence the above assumption 
is  equivalent to the Axiom of Choice.

In this note we will generally assume the Axiom of Choice.
 See Remark \ref{nochoi}  for some comments
about choiceless situations.
\end{definition}

\begin{problem} \labbel{undef}
Is it true that, for every $\Delta$-game $\mathbf G$ and 
 every choice of the First Player between Left and Right,
either First or Second has a winning strategy on $\mathbf G$?
Is this true at least for countable sums of combinatorial games?
 \end{problem}

\begin{fact} \labbel{fact}
Assume that $I$ is a set and  $\mathbf G = \dad _{i \in I} G_{i}$
is a $\Delta$-game.

If $I$ is finite, say, $I=\{ i_1, \dots, i_n  \}$, 
then options, runs, etc on $\mathbf G$ are 
defined exactly as in the classical theory of combinatorial games
on $G_{i_1} + \dots + G_{i_n}$.
\end{fact}

\begin{example} \labbel{dicot}
An infinite sum of impartial games 
is always second winner, provided infinitely many summands
are different from $0$.
Indeed, say in a countable sum,
no player can annihilate all the games in a finite number
of moves.
On the other hand, Second can turn to $0$ all the summands 
in $ \omega$ steps.  At the limit step First has to move, 
so she loses. In general, let Second well-order
the non $0$ summands with the order type of some cardinal $\kappa$.
Again, Second can turn to $0$ all the summands in $\kappa$ steps,
and then at the $\kappa$ limit step First cannot move.

More generally, by the same argument, an infinite sum of dicotic games 
is always second winner, provided infinitely many summands
are different from $0$.
Recall that a combinatorial game is \emph{dicotic}
if in every subposition different from $0$
every player has a move.  
 \end{example}

\begin{definition} \labbel{leftper}    
We call a combinatorial game $ G$  \emph{Left persistent}
if Left has an option in $ G$, and Left still has an option in
any subposition of $ G$ after any finite number of 
(consecutive) Right moves.
We are not asking that Left has an option after a Left move;
actually, $ G$ might even be a losing game for Left.
For example, if $G_0 = \{ \, -10 \mid  \,\}$,
$G_1 = \{ \, -10 \mid  G_0\,\}$, \dots,
$G_{n+1} = \{ \, -10 \mid  G_n \,\}$,
then each $G_i$ is a Left persistent combinatorial game. 
For every $X \subseteq \mathbb N$, also 
$G_ X  = \{ \, -10 \mid G_i  \, (i \in X) \,\}$  is a Left persistent combinatorial game.
All the above examples, except for $G_0$,  have negative value.

\emph{Right persistent} games are defined in the symmetrical way. 
\end{definition}

For example, every strictly positive number $G$ 
is Left persistent, since each Right move on $G$ increases
the value of $G$ \cite[Proposition 3.2]{S}. 
On the other hand, 
$\{ \, 0 \mid * \,\}$ (which is not a number) is strictly positive,
but not Left persistent. 
Moreover, $\{ \,  0\,\, ||\,\, 1|{-1} \,\}$ is
Conway equivalent to $1$, but is not
Left persistent. It follows that being  Conway equivalent 
to a positive number does not imply  Left persistence.

\begin{proposition} \labbel{lp} 
If  $\mathbf G= \dad _{n \in \mathbb N } G_{n}$
is a $\Delta$-game and  there are infinitely many indices
$i$ such that $G_i$ is Left persistent (resp., Right persistent), then Left
(resp., Right)
has a winning strategy in $\mathbf G$ when playing Second.  
\end{proposition}
  
\begin{proof}  
Left chooses a function $f: \mathbb N \to \mathcal P(\mathbb N)$
such that (a) $f(n)$ is infinite, for every $n \in \mathbb N$, (b)
if $i \in f(n)$, then  $G_i$ is Left persistent, and (c) $m \neq n$, 
$i \in f(n)$, $j \in f(m)$ imply $i \neq j$, namely, the $f(n)$s are 
pairwise disjoint.

For each Right move on some game $G_ \ell$,
Left answers with a move on $G_ i$,
for some $i \in f(\ell)$. 
Since each game $G_ \ell$ is a combinatorial game, 
a play on $G_ \ell$ ends after a finite number of moves, 
hence Left never ends up out of moves, since $f(\ell)$
is infinite and Left has a move
on each game in $f(\ell)$.
 It might happen that Right has already moved
on $G_i$, but  $G_i$ is Left persistent,
hence Left has indeed a move on it.

The point is that, under the above rules for alternating games,
Right always plays first at each limit stage, hence
each Left move immediately follows a Right move.
 \end{proof} 

So, for example,
$ \omega - \frac{1}{ \omega } +
\omega^ \omega  - \frac{1}{ \omega ^ \omega  } +
\omega^{ \omega^ \omega }  - \frac{1}{ \omega^{ \omega^ \omega } } + \dots$
 is a second winner game,
though each partial sum is a very large surreal!

If $|I|> | \mathbb N|$, one needs $|I|$-many
Left-persistent games in order to get the analogue of 
Proposition \ref{lp}. It seems that the axiom of choice
is needed, in this case. 

It follows from Proposition \ref{lp}
that the outcome of a $\Delta$-game is not
invariant under taking Conway equivalent summands.
For example, $-1+0+0+0+\dots$ is Right winner.
On the other hand, if $G$ is the Left Persistent 
and Right persistent game $\{ \, -1 \mid 1 \,\}$,
then  $-1+G+G+G+\dots$
is second winner, by 
Proposition \ref{lp}.
However, $G$ is Conway equivalent to $0$.

\section{Extended plays on finite sums of $\Delta$-games} \labbel{pfs}

As clear from the previous section, 
in case  there are infinite runs,
starting some infinite game as
a First player is a sort of handicap, if ``First always
remains First at limit steps''. 
Suppose that $\mathbf G$ and $\mathbf H$
are $\Delta$-games, $\mathbf G$ is second winner and,
say,
 Left can win on $\mathbf H$ 
when playing first.
If $\mathbf G$ and $\mathbf H$ are
classical combinatorial games, Left has a winning strategy
by first playing on $\mathbf H$ and then using her winning strategy
on $\mathbf H$ as a First player each time Right replies moving on $\mathbf H$,
and
using her winning strategy
on $\mathbf G$ as a Second player each time Right replies moving on $\mathbf G$.
In the case of infinite games, a tentative strategy for Left
would still be to play First on $\mathbf H$. If  Right always plays
on $\mathbf G$ at his turn on the first $ \omega$ moves,  
then Left is on her turn at the $ \omega$th move. But Right 
has never moved on $\mathbf H$, hence it is not necessarily the case
that Left has still a winning move on the remainder of $\mathbf H$ (this would
be her second consecutive move on $\mathbf H$; we are just
assuming that Left wins playing first on $\mathbf H$).
Left has not necessarily a winning move on $\mathbf G$, either,
since, we assume that $\mathbf G$ is second winner, so, by our
conventions, Left has a winning strategy on $\mathbf G$ if she is always 
 Second after a limit step.

The above example itself suggests  a possible solution.
In the example, only one move has been performed on $\mathbf H$
and infinitely moves have been performed on $\mathbf G$.
Hence it is natural to ask that, after $ \omega$ moves,
the player to move is the first player who has made some move 
on a game \emph{on which infinitely moves have been performed}.
For limit ordinals larger than $ \omega$, we consider
``cofinally many moves'' on some game, rather than infinitely many
moves.
In this way we can reconstruct many salient aspects of Conway theory! 

Notice that, if we apply the above ideas, the notion of a run
is not modified: the only notion which is changed is the notion of
an alternating run. Now for the explicit details.

\begin{definition} \labbel{oplg}   
Recall the assumptions from item (B) in Definition \ref{sgam}.
We have a sequence  $( \mathbf G_j) _{j \in J} $  of $\Delta$-games,
with   $\mathbf G_j = \dad _{i \in I_j} G_i$, for each $j \in J$ and we assume
 that the sets $I_j$ are pairwise disjoint.
We have set   $\daf _{j \in J}  \mathbf G_j = \dad _{i \in I} G_i$,
where $I = \bigcup _{j \in J} I_j $.  

(a)
Assume further that $J$ is finite, say,
$J= \{ 0,1, \dots, n-1 \} $.
Denote by $\mathbf G_0 \oplus \mathbf G_1 \oplus
\dots \oplus \mathbf G_{n-1}  $ the sum  
$\daf _{j \in J}  \mathbf G_j = \dad _{i \in I} G_i$,
indexed by $I=\bigcup _{j \in J} I_j $.
So far, $\circled{ G} = \mathbf G_0 \oplus 
\dots \oplus \mathbf G_{n-1}  $
 is just an ordinary $\Delta$-game.
We endow  $\circled{ G}$ with a partition 
of $I$ onto finitely many classes:
the classes are exactly $I_0, \dots, I_{n-1}$. 

An \emph{extended game}, or   \emph{$\oplus$-game}, is a $\Delta$-game
 indexed on some set $I$, together with
a partition 
of $I$ onto finitely many classes.
Every finite sum of $\Delta$-games becomes
naturally a $\oplus$-game, as in the previous paragraph.
Conversely, every $\oplus$-game can be written
as a finite $\oplus$-sum as above, just group together
in a single $\Delta$-game all the combinatorial games
corresponding to any given class of the partition. 
Rather than indicating explicitly the partition, we will denote
$\oplus$-games as above, as $\mathbf G_0 \oplus 
\dots \oplus \mathbf G_{n-1} $. 
A natural (finitary) associative operation $\oplus$
is defined on the class of $\oplus$-games.

(b)
\emph{Moves, runs, resultants} for $\oplus$-games
are defined  as in the case of $\Delta$-games; 
in such cases  the partition
 plays no role.

(c)
However, we will use a different notion
of an alternating run. An \emph{extended alternating run}
$r$ on a $\oplus$-game $\circled{ G}$ is defined as follows.

Players play alternatively on successor moves, as in ordinary
games and in $\Delta$-games; in particular, there is a First player
which starts moving on $\circled{ G}$, then the Second
player makes the second move, subsequently First moves
again and so on. 

At any limit step,
say, for some limit ordinal $\beta$, we define some 
partition class $I_j$ to be \emph{active} 
(with respect to $\beta$ and to some run $r$ defined up to $\beta$
on $\circled{G}$)
if, for every $\gamma < \beta $,
there is   some $\gamma'$ such that  
$\gamma \leq \gamma ' <  \beta $
and the move $r( \gamma ')$ has been made on 
some $G_i$ with   $i \in I_j$. Notice that, since $\beta$ 
is limit and there is only a finite number of classes
in the partition, there is at least one active partition.

Then, for $r$ to be considered extended alternating, 
we require that the move $r ( \beta )$ is done by
the first player who has ever moved on a game
whose index lies on an active partition. Namely, if $ \delta $
is the first ordinal such that $r ( \delta )$ is a move in a game
whose index lies on an active partition, then we require that 
the Same player which has made the move   
$r ( \delta )$ makes the move $r ( \beta  )$.
Then the players alternate until the next limit step, if
any, thus the Other player makes the move 
$r ( \beta + 1 )$, the Same player makes the move 
$r ( \beta +2 )$ and so on.

(d)
As usual, a \emph{play} (or an \emph{extended play},
when we want to make clear that we are applying the above convention) 
is a maximal extended alternating run.

For example, if $\circled{G}= \mathbf G_1 \oplus \mathbf G_2$,
Right moves  first on $\mathbf G_1 $,
Left moves  second on $\mathbf G_2 $  and subsequently all players
move   on $\mathbf G_2 $ (for the first
$ \omega$ moves), then, according to the extended rule,
Left is to move at the $ \omega$th step. 

(e)
Notice that the above notion of an extended play does not
necessarily respect the
$\oplus$ sum. 
Suppose that 
$\circled{G}= \mathbf G_1 \oplus \mathbf G_2$,
$\circled{H}= \mathbf H_1 \oplus \mathbf H_2$
$\circled{K}= \mathbf K_1 $
and the play is made on 
$\circled{G} \oplus \circled{H} \oplus \circled{K} =
 \mathbf G_1 \oplus \mathbf G_2 \oplus
 \mathbf H_1 \oplus \mathbf H_2 \oplus \mathbf K_1 $.
Suppose that Right moves first on $\mathbf K_1$,
Left moves second on $\mathbf G_1$,
Right moves third on $\mathbf H_1$,
and then the players alternatively move 
on $\mathbf G_2$ and $\mathbf H_1$.    

In the sum 
$\circled{G} \oplus \circled{H} \oplus \circled{K}$
the ``active'' $\oplus$-games are 
$\circled{G} $ and $  \circled{H} $    
and Left is the first to have made a move in one
of those two games.
But the partition given by 
Definition \ref{oplg} has actually $5$  
classes, the active partitions in
$\mathbf G_1 \oplus \mathbf G_2 \oplus
 \mathbf H_1 \oplus \mathbf H_2 \oplus \mathbf K_1 $
correspond to the games $\mathbf G_2$ and $\mathbf H_1$.
The first player having made a move on
one of those games is Right, hence, according
to the extended rule, Right is to move at the 
$ \omega$th step.  
 \end{definition}

We need a further definition to get
an analogue of Conway theory to work.

\begin{definition} \labbel{unif} 
We say that \emph{Left has a uniform winning strategy
playing as Second on some $\oplus$-game $\circled{G}$}
if Left, playing Second, has a winning strategy on  $\circled{G}$
with the further requirement that Second is never the first to move
at limit stages, namely, if Left applies such a strategy, then
extended runs turn out to be ``ordinary'' runs 
in the sense of Definition \ref{run}. 

The definition of a \emph{Right uniform winning strategy}
is symmetrical. An $\oplus$-game $\circled{G}$ is a
\emph{uniformly second winner game} if  
Left have a uniform winning strategy
when playing as Second on $\circled{G}$ and 
 Right, as well,   has a uniform winning strategy
when he plays as Second.
\end{definition}   

\begin{problem} \labbel{mainpro}
Are there $\oplus$-games such that Left,
when playing Second, has a winning strategy but
not a uniform winning strategy?

In particular, are there  second winner $\oplus$-games
which are not
uniformly second winner game?
 \end{problem}  

The property in Definition \ref{unif} is used in order
to show that the sum of two  
uniformly second winner games is still second winner.
Indeed, 
in the case of alternating runs defined as in Definition \ref{run},
it is elementary to see that the sum (in the sense of 
Definition \ref{sgam}(B)) of any family of second winner
$ \Delta $-games is still second winner. 
However, we do not know how to prove this in the case
of extended plays on $\oplus$-games since, say
in a sum $\circled{G} \oplus \circled{H}$,
it might happen that the winning strategies for Second
in these games, considered as separate games,
lead her to play first at some limit stage in both games.
But this is impossible when she plays
in $\circled{G} \oplus \circled{H}$,
 hence First might  play first even at this stage
on one of the two games, and this might ruin Second's 
strategy in that game.
The assumption of uniformity from Definition \ref{unif}
prevents the above situation to happen.

\begin{lemma} \labbel{lemmon}
  \begin{enumerate}    
\item 
If Left has a uniform winning strategy
playing as Second on each of the $\oplus$-games $\circled{G}$
and $\circled{H}$, then Left has a uniform winning strategy
playing as Second on $\circled{G} \oplus \circled{H}$.
\item 
If the $\oplus$-games $\circled{G}$
and $\circled{H}$ are both
uniformly second winner $\oplus$-games,
then their sum $\circled{G} \oplus \circled{H}$
is a uniformly second winner $\oplus$-game.
  \end{enumerate} 
 \end{lemma}

 \begin{proof} 
(1) Whenever Right plays on the $\circled{G}$ component,
Left replies using her uniform winning strategy on
the $\circled{G}$ component,
and symmetrically for $\circled{H}$. Since both strategies
are uniform, 
Right is always the first to move at each limit stage
in the play on $\circled{G} \oplus \circled{H}$,
hence Left can go on applying her strategies.

(2) is immediate from (1) and the symmetric version dealing with Right.
\end{proof}

\begin{problem} \labbel{probboh}
Are the following itemized statements true?
Of course, as usual when dealing with
$\oplus$-games, we always assume extended play. 

  \begin{enumerate}[(i)] 
\item 
Is the outcome of every $\oplus$-game
always defined?    
  \end{enumerate} 
Let  $\circled{G}$
and $\circled{H}$ be two $\oplus$-games and assume that
Left has a uniform winning strategy
playing as Second on $\circled{G}$.

 \begin{enumerate}[(i)]   
\setcounter{enumi}{1}
\item
 Left has a  winning strategy
playing as Second on $\circled{G} \oplus \circled{H}$
if and only if 
Left has a  winning strategy
playing as Second on $\circled{H}$.
\item
 Left has a  winning strategy
playing as First on $\circled{G} \oplus \circled{H}$
if and only if 
Left has a  winning strategy
playing as First on $\circled{H}$.
\item
The outcome of $\circled{G} \oplus \circled{H}$
is defined if and only if 
the outcome of  $\circled{H}$
is defined and, if both are defined,
they are equal.
   \end{enumerate} 

As a standard argument, 
if (ii)-(iii) are true, then (iv) is true.
 \end{problem}

\section{Regular alternating plays} \labbel{reg}

\begin{definition} \labbel{regg}
(a) Suppose that $G$ is a combinatorial game
and $s$ is a run on $G$.   
A \emph{final subrun} $f$ of $s$
is the sequence of the last $m$ moves of $s$,
for some $m$.  Formally, $f$ is a run
on the position $s(n-m)$, where $n$ is the length of $s$.      
The (possibly empty) \emph{maximal alternating final subrun}, for short,
\emph{maf subrun} of some run $s$ is a final alternating subrun   $f$
which has even length and is maximal with such properties.

For example, if the players move
in the following order on $s$:
LRLRRLLRRLRLRL,
then the maf subrun is given by the last $6$ moves.
If the order is 
  LRLRLRRLLRLRL
the maf subrun is given by the last $4$ moves,
since we require its length to be even.
Finally, 
the maf subrun
of LRLRLRLRLRLL
 is the empty run of length $0$. 

(b) 
We define the regular alternating play 
on some $\Delta$-game $\mathbf G$ 
by not considering maf subruns 
(induced on summands)
when determining which player should move at limit stages.
In more detail:

(c) A run $r$ on a $\Delta$-game $\dad _{i \in I} G_{i}$  
naturally  induces a run $r_i$ on each summand $G_i$,
simply considering the subsequence of $r$
consisting of those moves which are played on $G_i$.    

(d) A move of $r$ is \emph{regular}  if the corresponding move  
in $r_i$, as defined above (notice that $i$ is uniquely determined),
does not belong to  a  maf subrun of $r_i$.
 
(e)  The \emph{regular alternating play} 
on a $\Delta$-game $\dad _{i \in I} G_{i}$
is played as follows.
As usual, the two players move alternately
at successor stages.
At limit stages, one considers the sequence
of the already played regular moves.
If this sequence has not a last element, that is, it
is indexed by a limit ordinal, then First has to move.
If the sequence has a last element and Some player has moved 
at this last element, then the Other player has to move.
  \end{definition}   

The notion of a regular alternating play alleviates some paradoxical aspects
of $ \Delta $-games. 
It follows from Proposition \ref{lp} that,
assuming the alternating play from Definition \ref{run}(c),
a countable sum which contains infinitely many copies of
the Conway $0$-game $\{ \, -1 \mid 1 \,\}$ is second winner, no matter
the remaining games. 
On the other hand, 
$1$ together with any number of copies of 
$\{ \, -1 \mid 1 \,\}$ is Left winner under the regular play.
Suppose, say, that Left plays first. She can play on 
$1$ and then reply to any Right move on the same game.
Under the regular rule, such pairs of moves are not ``counted''
at limit stages, so that Right is always the first to move at any limit stage. 
 More generally, we are going to prove in  Proposition \ref{zeronon},
that, under the regular alternating play,
second winner, that is, Conway $0$
  combinatorial games do not affect the outcome of a
$\Delta$-game. 

As another example, and in contrast with Example \ref{dicot},
under the regular play, any nonempty sum whose summands
are all $\uparrow = \{ \, 0 \mid * \,\}$ is Left winner.
If Left plays as second, let her always reply on 
the game on which Right has just played.
If Left plays as first, let her move on some $ \uparrow $
and then play as above. Because of the  
  regular rule, Right has to move first at each limit step,
so he eventually loses.

\begin{proposition} \labbel{zeronon}
Under the regular alternating  play, Conway zero combinatorial games
do not affect the outcome of a $\Delta$-game. 

In detail, under the regular alternating  play,
if $\dad _{i \in I} G_{i}$ is a $\Delta$-game,
 $J \subseteq I$ is a set  such that  
$G_i$ is a second winner combinatorial game,
for every  $i \in J$,
then  the outcome of $\dad _{i \in I} G_{i}$
is defined if and only if 
the outcome of $\dad _{i \in I \setminus J} G_{i}$
is defined and, if this is the case, the two outcomes
are the same.
 \end{proposition}  

\begin{proof}
If some player $P$ has a winning strategy on 
$\dad _{i \in I \setminus J} G_{i}$, then $P$ has
a winning strategy on 
$\dad _{i \in I} G_{i}$ by using the strategy 
for $\dad _{i \in I \setminus J} G_{i}$ when playing
on games in $I \setminus J$, and using the second winner strategy 
on each $G_i$, for $i \in J$. Namely, $P$ always plays
on games in  $I \setminus J$, except immediately after 
an adversary move on some game in $J$.

Conversely, assume that $P$ has  a winning strategy on 
$\dad _{i \in I} G_{i}$. When the strategy involves
$P$ moving  as a first player on some game $G_i$ with  $i \in J$,
let the adversary reply in $G_i$ using the second winner strategy
for $G_i$. Forgetting about the moves on games in $J$,
we get a winning strategy for $P$ in $\dad _{i \in I \setminus J} G_{i}$. 
 In fact, the proof shows that if $P$ 
has  a winning strategy on 
$\dad _{i \in I} G_{i}$, then $P$ has a winning strategy in which
she never plays  first on games in $J$.
 \end{proof}

As shown by the next proposition,
paradoxical aspects remain also with the regular alternating play.

\begin{proposition} \labbel{pnumb}
If  $\mathbf G= \dad _{n \in \mathbb N } G_{n}$
is a $\Delta$-game and  there are infinitely many indices
$i$ such that $G_i$ is a strictly positive number 
(resp., a strictly negative number), then Left
(resp., Right)
has a winning strategy in $\mathbf G$ when playing Second,
with the regular alternating play.  
 \end{proposition} 

\begin{proof} 
Similar to the proof of Proposition \ref{lp}.
The First player could try to apply the strategy of always
replying on the game the Second player has moved on,
so that Second could be forced to be the first to play at some limit step.
But Second can always avoid this by temporarily
pausing her main strategy and playing on the same positive
number until First has no move on it.
\end{proof}   

Hence adopting the regular alternating play is only palliative.
For example, by Proposition \ref{pnumb},
the same $\Delta$-game mentioned after
Proposition \ref{lp}, that is,
$ \omega - \frac{1}{ \omega } +
\omega^ \omega  - \frac{1}{ \omega ^ \omega  } +
\omega^{ \omega^ \omega }  - \frac{1}{ \omega^{ \omega^ \omega } } + \dots$
 is a second winner game
also under the regular alternating play.

It might happen that the  arguments in Proposition \ref{zeronon}, together 
with a strategy stealing argument, can be used to prove the following
result, but we have not yet fully checked details.

\begin{conjecture} \labbel{regthm}   
Under the regular alternating  play,  the outcome of a $\Delta$-game
is invariant with respect to Conway equivalence of the summands. 

In detail,
suppose that $( G_i) _{i \in I} $ 
and $( H_i) _{i \in I} $
are combinatorial games and
$G_i$ is Conway equivalent to $H_i$,  
for every $i \in I$.
Then, 
 under the regular alternating  play,
 the outcome of $\dad _{i \in I} G_{i}$
is defined if and only if 
the outcome of $\dad _{i \in I } H_{i}$
is defined and, if this is the case, the two outcomes
are the same.
\end{conjecture}

\section{Further remarks} \labbel{fur} 

\begin{problem} \labbel{comb}
Can we combine the ideas from 
Sections \ref{pfs} and \ref{reg}.
Namely, is there some useful notion of an
``extended regular'' alternating play?

Another possibility which can be considered 
(``Annihilating zeros'') has been hinted in
\cite[Section 5.4]{scg}.   
 \end{problem}

\begin{remark} \labbel{nochoi}
(a) It seems that that the Axiom of Choice plays a relevant role
in all the above arguments.
In choiceless situations, there might exist an
infinite set (= not equipotent with a finite set) $D$ 
which is not Dedekind infinite, that is, there is no injection
from $ \mathbb N$ to $D$.

If $D$ is such a set and 
$\mathbf G =\dad _{i \in D} G_{i}$
is a $\Delta$-game,
then there is no infinite run on  
$\mathbf G$, since this would provide
an injective function from $ \mathbb N$ to $D$.

However,  we can still talk of winning strategies.
For example, if each $G_i$ 
is $\{ \, * \mid * \,\}$, then Second has a winning strategy
by always moving on the game on which First has just moved.   

If each $G_i$ 
is $\{ \, * \mid * \,\}$, except that just one
$G_{\bar i}$ is $*$, then First has a winning strategy by 
making the first move on  $G_{\bar i}$ 
and then replying on the game Second moves on.
(This contrasts with the case when $I$ can be put in a 
bijection with $ \omega$. In this case, $\mathbf G =\dad _{i \in D} G_{i}$
is second winner, by Example \ref{dicot}.)

(b) Note that some care is needed in defining a winning strategy.
Strategies here cannot be total functions, since otherwise, as above, 
we could get an injective function from $ \mathbb N$ to $D$.

Hence we define a \emph{strategy}  to be a partial function
from the set of subpositions of  $\mathbf G$.
A strategy $\sigma$ is a \emph{winning strategy}, say, for Left playing
as Second, if, for each even ordinal $\alpha$ and every alternating run
of length $\alpha-1$ with resulting
subposition $\mathbf H$, $\sigma$ furnishes     
a legal move for Left on $\mathbf H$,
\emph{provided that, up to $\alpha-1$, Left has played the run 
according to $\sigma$}.  
 \end{remark}

\end{document}